\documentclass{article}
\usepackage[francais]{babel}
\usepackage[latin1]{inputenc}
\usepackage[T1]{fontenc}
\usepackage{amsmath}
\usepackage[psamsfonts]{amsfonts}
\usepackage{amsthm}
\usepackage{amsopn}
\usepackage{amscd}
\usepackage{makeidx}
\usepackage[dvipdf]{graphics}
\theoremstyle{plain}
\newtheorem{Theo}{Théorème}
\theoremstyle{remark}

\newtheorem*{Demo}{Démonstration de la propriété 2}
\newtheorem*{Demoo}{Démonstration de la propriété 3}

\theoremstyle{definition}

\catcode`\«=\active \catcode`\»=\active
\def«{\og\ignorespaces}
\def»{{\fg}}
\normalsize

\begin{document}

\begin{center}
{\LARGE\bf Inégalité d'Ahlfors en dimension supérieure}\footnote{Published in Mathematische Annalen:  \\
http://www.springerlink.com/content/u870536335x7468w/?p=f03b483c6abc415fa2f1f22220607a70\&pi=0 \\
The original publication is available at www.springerlink.com .}
\end{center}

\begin{center}
{\large Benoît Saleur}\footnote{Département de 
Mathématiques de la faculté des sciences d'Orsay, Université Paris-Sud 11, 91405 Orsay Cedex. \\ 
Adresse électronique: saleur@clipper.ens.fr \\
Mots clefs: Hyperbolicité complexe, inégalité isopérimétrique, surface de Riemann.\\
Clefs AMS: 30F45, 32Q45.
}
\end{center}

\normalsize \vspace*{0.5cm} \hrule\vspace*{1.5cm}

\begin{center}
{\bf Abstract}
\end{center}

\begin{center}
{ We prove an Ahlfors' like inequality for the holomorphic curves with boundary of a complex compact Kobayashi-hyperbolic manifold.}
\end{center}

\section{Introduction}

L'inégalité d'Ahlfors (voir par exemple \cite{DeThélin}, \cite{Nevanlinna}) contient sous une forme géométrique l'essentiel de la théorie de distribution des valeurs des fonctions analytiques. Dans le cas d'une surface d'arrivée sans bord, elle s'énonce ainsi:

\begin{Theo}
Soit $\Sigma_{0}$ une surface de Riemann compacte munie d'une métrique hermitienne fixée. Alors il existe une constante $h>0$ telle que pour toute surface de Riemann compacte à bord $\Sigma$ et toute fonction $\Sigma \stackrel{f}{\rightarrow} \Sigma_{0}$ holomorphe à l'intérieur de $\Sigma$ et lisse au bord, l'inégalité suivante soit vérifiée:
$$\min (0,\chi (\Sigma)) \leq \chi (\Sigma_{0}) S + hL.$$
On a noté $S =\frac{\operatorname{Aire}(f(\Sigma ))}{\operatorname{Aire}(\Sigma_{0} )}$ le nombre moyen de feuillets au-dessus de $\Sigma_{0}$, $L=\operatorname{Longueur}(f(\partial \Sigma ))$ la longueur du bord  de $f(\Sigma )$, et $\chi (\Sigma_{0})$, $\chi (\Sigma )$ les caractéristiques d'Euler-Poincaré des surfaces $\Sigma_{0}$ et $\Sigma$.
\end{Theo}

On se propose ici d'établir une inégalité de type Ahlfors pour les courbes holomorphes d'une variété complexe compacte hyperbolique au sens de Kobayashi, c'est à dire ne contenant pas de courbe entière. On fixera donc une variété complexe compacte hyperbolique $X$ munie d'une forme hermitienne $\omega$. \\
Une inégalité isopérimétrique linéaire pour les disques holomorphes a déjà été établie par Julien Duval (voir \cite{Duval}), à l'aide d'un lemme à la Brody. Elle est également une conséquence implicite d'un résultat établi par Bruce Kleiner  dans son preprint "Hyperbolicity using minimal surfaces"(voir \cite{Kleiner}). Le lemme de Brody a également permis à Jean-Pierre Demailly d'obtenir une inégalité de type Ahlfors pour les courbes holomorphes sans bord (voir \cite{Demailly}). L'inégalité ici démontrée combine ces deux résultats. La démonstration repose sur des estimations longueur-aire et sur le théorème de Gauss-Bonnet, et emprunte beaucoup aux méthodes utilisées par Bruce Kleiner (voir \cite{Kleiner}).\\

Soient $\Sigma$ une surface de Riemann compacte à bord, de caractéristique d'Euler-Poincaré $\chi$, et $f: \Sigma \mapsto X$ une fonction holomorphe sur l'intérieur de $\Sigma$ et lisse au bord (on parlera de courbe holomorphe à bord). On munit $\Sigma$ de sa métrique de Poincaré $h$, de courbure $-1$. La pseudo-métrique $g=f^{*}\omega $ est conforme à $h$: il existe sur $\Sigma$ une fonction lisse $\Phi$ telle que: $g =\Phi^{2}h$. De plus, par le lemme de Brody (voir \cite{Brody}) il existe une constante $M\geq \sqrt{2}$ ne dépendant que de $(X, \omega )$ telle que $\Phi \leq \frac{M}{\sqrt{2}}$.

Nous pouvons alors énoncer le théorème suivant:

\begin{Theo}
L'une des deux inégalités suivantes est vérifiée:
$$\operatorname{Aire}_{\omega}(f(\Sigma ))\leq -2\pi M^{2} \min (0,\chi (\Sigma))$$
ou
$$\operatorname{Aire}_{\omega}(f(\Sigma ))\leq 2e^{8M^{2}} \operatorname{Longueur}_{\omega}(f(\partial \Sigma)).$$
\end{Theo}

Ce résultat peut s'interprêter comme une notion faible de courbure négative pour la métrique de Kobayashi.

\section{Démonstration du théorème}

La démonstration s'inspire de la méthode utilisée par Fiala pour prouver l'inégalité isopérimétrique dite de Bol-Fiala (voir \cite{Fiala}).\\

Commençons par remarquer que la courbure de Gauss $K_{g}$ de la pseudo-métrique $g=f^{*}\omega$, définie hors de l'ensemble critique de $f$, peut être supposée inférieure ou égale à $1$: en effet, $f(\Sigma )$ est alors localement une sous-variété complexe de $X$, et sa courbure est donc majorée par une constante (voir l'article de J. Lafontaine dans \cite{Audin} pour plus de détails). Il suffit alors de normaliser $\omega$ pour obtenir $K_{g}\leq 1$.\\

Il sera nécessaire par la suite de disposer d'une vraie métrique: en posant $g_{\varepsilon} = (\Phi + \varepsilon )^{2}h$, on obtient une famille de métriques très proches de $g$, dont la courbure vérifie: 
$\overline{\displaystyle{\lim_{\varepsilon \rightarrow 0}}}K_{g_{\varepsilon}} \leq 1$. Il suffit alors de démontrer le théorème pour les $g_{\varepsilon}$ et de passer à la limite. Par souci de clarté, nous supposerons directement que $g$ est une métrique.\\

Notons $d$ la distance associée à la métrique $g$. Pour tout réel $t>0$, posons: $$\Sigma (t) = \{ p \in \Sigma \text{ | } d(p, \partial \Sigma ) \geq t\}.$$
Les frontières $\partial \Sigma (t)$ de ces domaines sont appelées "courbes parallèles"' au bord $\partial \Sigma$. Elles peuvent présenter des points angulaires: ceux-ci apparaissent à un temps $t_{0}$ et se propagent aux temps $t>t_{0}$. Deux phénomènes causent l'apparition de points angulaires au temps $t_{0}$ (voir \cite{Arnold} ainsi que l'article de M.P. Muller dans \cite{Audin}):
\begin{enumerate}
  \item L'existence de points de contact pour la courbe $\partial \Sigma (t_{0})$, c'est à dire de points $P$ de $\Sigma$ pour lesquels il existe deux points distincts $P_{1}$ et $P_{2}$ de $\partial \Sigma$ et deux arcs de géodésiques de même longueur $t_{0}$, orthogonaux à $\partial \Sigma$, et reliant $P$ à $P_{1}$ et $P_{2}$. Pour tout $t>t_{0}$, la courbe $\partial \Sigma (t)$ présente un point angulaire issu du point de contact. La topologie de $\Sigma (t)$ est modifiée au passage de telles valeurs de $t$. 
  \item L'existence de points en lesquels la courbure géodésique est infinie, c'est à dire de points singuliers pour la courbe $\partial \Sigma (t_{0})$. Si $P$ est un point singulier de $\partial \Sigma (t_{0})$, la courbe $\partial \Sigma (t)$ présente un point angulaire issu de $P$ pour tout $t>t_{0}$.
\end{enumerate}

En un point angulaire d'une parallèle $\partial \Sigma (t)$, l'angle entre les normales entrantes est inférieur à $\pi$. En effet, pour un réel $\varepsilon >0$ assez petit, la courbe $\partial \Sigma (t)$ est l'enveloppe des frontières des boules de rayon $\varepsilon$ et de centres les points de $\partial \Sigma (t-\varepsilon )$, et en un point d'intersection entre deux boules, l'angle entre les normales sortantes est inférieur à $\pi$.\\
D'après \cite{Arnold}, on peut approximer la courbe $\partial \Sigma$ par une courbe lisse en position générale pour les points doubles des parallèles comme pour les points singuliers. On supposera donc par la suite qu'il n'y a qu'un nombre fini de points angulaires.\\
Notons $\chi (t)$ la caractéristique d'Euler-Poincaré de la surface $\Sigma (t)$. On constate que la fonction $\chi$ est croissante. En effet, il existe deux types de contacts:

\begin{enumerate}
  \item Un contact entre deux composantes connexes distinctes de $\partial \Sigma (t)$. Alors $\Sigma (t^{+})$ a une composante de bord de moins que $\Sigma (t^{-})$.
  \item Un contact d'une même composante connexe de $\partial \Sigma (t)$ avec elle-même. Alors on obtient la surface $\Sigma (t^{+})$ en retirant un disque à deux trous de $\Sigma (t^{-})$.
\end{enumerate}
Dans les deux cas, on a bien $\chi (t^{+})= \chi (t^{-}) +1$.\\

Par la suite on notera $a(t)$ l'aire de la surface $\Sigma (t)$ et $l(t)$ la longueur de son bord $\partial \Sigma (t)$ (pour la métrique $g$). La fonction $a$ est dérivable, et on a la relation classique: $a'(t)=-l(t)$. La fonction $l$ est continue (car la courbe $\partial \Sigma$ est en position générique) et dérivable à droite (la notation $l'$ désignera la dérivée à droite de $l$).\\
On notera enfin $\chi^{-} = \min (0, \chi )$.\\

Les trois propriétés suivantes sont alors vérifiées et permettent de conclure:
\begin{enumerate}
  \item $a'=-l$.
  \item $a-l'\geq 2\pi \chi^{-}$.
  \item $\forall T\geq 0$, $a(T) + \pi M^{2}\chi^{-} \leq \frac{M^{2}}{2\int_{0}^{T}\frac{dt}{l(t)}}$.
\end{enumerate}

Commençons par établir le Théorème $2$ avant de démontrer les propriétés $2$ et $3$:\\

Supposons que $a(0)\geq 2 \pi M^{2} \chi^{-}$. 
En intégrant la relation $l' \leq a - 2\pi \chi^{-}$ entre $0$ et $T$ et en majorant $a(t)$ par $a(0)$ pour tout $t\geq 0$, on obtient: $$l(T) \leq l(0) + T (a(0) - 2\pi \chi^{-})
\leq l(0) + 2Ta(0).$$
D'une part, on en déduit:
$$\frac{1}{\int_{0}^{T}\frac{dt}{l(t)}} \leq \frac{2 a(0)}{\log (1 +  2T \frac{a(0)}{l(0)})}.$$
D'autre part, en intégrant la relation $1$, on obtient l'inégalité suivante, valable même lorsque $a(T)=0$ et $l(T)=0$:  $$a(T)\geq a(0)-Tl(0) -T^{2}a(0).$$
L'inégalité $3$ s'écrit alors, en notant $x= \frac{a(0)}{l(0)}$:
$$x-T-T^{2}x + \frac{\pi M^{2} \chi^{-}}{l(0)} \leq \frac{M^{2} x}{\log (1  + 2T x)}.$$
Comme $a(0) \geq -2\pi M^{2}\chi^{-}$, cela s'écrit:
$$ x(\frac{1}{2}-T^{2})-T  \leq \frac{M^{2}x}{\log (1  + 2T x)}.$$ 
Posons $T=\frac{1}{2}$. Distinguons deux cas: soit le terme de gauche de l'inégalité ci-dessus est négatif ou nul, auquel cas $x\leq 2$, soit il est strictement positif, auquel cas:
$$\log \left (1 + x \right ) \leq \frac{M^{2}x}{x(\frac{1}{2}-T^{2}) - T}\leq \frac{2M^{2}x}{x-2}=2M^{2} + \frac{4M^{2}}{x-2}.$$
Pour $x\geq 4$, cette inégalité devient:
$$x\leq 2e^{8M^{2}}$$
comme annoncé.\\

\begin{Demo}
Elle est une conséquence presque immédiate du théorème de Gauss-Bonnet et suit l'article de M.P. Muller dans \cite{Audin}.
On note $\alpha_{i}$ les valeurs des éventuels angles intérieurs du bord. Le théorème de Gauss-Bonnet s'écrit alors:
$$\int_{\Sigma (t)}K_{g}v_{g} + \int_{\partial \Sigma (t)}k ds_{g} +\sum_{i} (\pi - \alpha_{i}) = 2\pi \chi (t)$$
où $K_{g}$ est la courbure de la métrique $g$ et $k$ est la courbure géodésique de la courbe $\partial \Sigma (t)$.\\
Or, un dessin montre que $l$ est dérivable à droite et que: $\int_{\partial \Sigma (t)}k ds_{g} = -l'(t) -2\sum_{i} \operatorname{cotan} (\frac{\alpha_{i}}{2})$. Comme $K_{g}\leq 1$, on a:
$$a(t) -l'(t) \geq 2\pi \chi (t) +2\sum_{i} \left (\operatorname{cotan}\frac{\alpha_{i}}{2}-\frac{\pi -\alpha_{i}}{2}\right ).$$
0r, d'une part, $0\leq \alpha_{i} <\pi$, soit $\displaystyle{\sum_{i}} \left (\operatorname{cotan}\frac{\alpha_{i}}{2}-\frac{\pi -\alpha_{i}}{2}\right ) \geq 0$, et d'autre part $\chi (t) \geq \chi^{-}$. La propriété est donc démontrée.
\end{Demo}

\begin{Demoo}
Elle découle de l'existence d'une inégalité isopérimétrique pour la métrique de Poincaré $h$ de courbure $-1$ et suit celle de B. Kleiner dans \cite{Kleiner}.\\
Notons $v_{g}$ la forme d'aire de la métrique $g$. Notons $a_{h}(t)$ et $l_{h}(t)$ respectivement l'aire de $\Sigma (t)$ et la longueur de $\partial \Sigma (t)$ pour la métrique hyperbolique $h$. Soit $t>0$. Comme $a_{h}(t)= \int_{\Sigma (t)}\Phi^{-2}v_{g}$, alors $a_{h}'(t)=-\int_{\partial \Sigma (t)}\Phi^{-2} ds_{g}$. l'inégalité de Cauchy-Schwarz s'écrit:
$$l_{h}(t)^2 = \left (\int_{\partial \Sigma (t)}\Phi^{-1}ds_{g}\right )^2 \leq \left (\int_{\partial \Sigma (t)} \Phi^{-2} ds_{g} \right )\left ( \int_{\partial \Sigma (t)}ds_{g} \right ) = -a_{h}'(t) l(t)$$
soit encore: 
$$\frac{-a_{h}'(t)}{l_{h}(t)^{2}}\geq \frac{1}{l(t)}.$$

On fixe un réel $T\geq 0$. Comme $a_{h}$ est décroissante, soit $a_{h}(T) + 2\pi \chi^{-} \leq 0$, auquel cas l'inégalité est triviale pour $T$, soit pour tout $0\leq t \leq T$, $a_{h}(t)+2\pi \min \chi^{-} >0$. L'inégalité isopérimétrique vérifiée par la métrique de Poincaré s'écrit: $a_{h}(t) \leq l_{h}(t) -2\pi \chi^{-}$ (voir par exemple \cite{Burago}). En l'injectant, on obtient:
$$\frac{1}{l(t)} \leq \frac{-a_{h}'(t)}{l_{h}(t)^{2}} \leq \frac{-a_{h}'(t)}{(a_{h}(t) + 2\pi \chi^{-})^2},$$
 ce qui donne par intégration:
$$\int_{0}^{T}\frac{dt}{l(t)} \leq \int_{0}^{T} \frac{-a_{h}'(t)}{(a_{h}(t) + 2\pi \chi^{-})^2}dt \leq \frac{1}{a_{h}(T) +2\pi \chi^{-}} - \frac{1}{a_{h}(0) + 2\pi \chi^{-}}\leq \frac{1}{a_{h}(T) +2\pi \chi^{-}}.$$
La propriété suit.
\end{Demoo}

\end{document}